\documentstyle{amsppt}
\magnification=1200
\hsize=150truemm
\vsize=224.4truemm
\hoffset=4.8truemm
\voffset=12truemm

\NoRunningHeads
 
\define\R{{\bold R}} 
\define\C{{\bold C}}
\let\th\proclaim
\let\fth\endproclaim

\newcount\refno
\global\refno=0

\def\nextref#1{
      \global\advance\refno by 1
      \xdef#1{\the\refno}}

\def\bref {\ref\global\advance\refno by 1\key{\the\refno}}

\nextref\BOT
\nextref\DUV
\nextref\GRO
\nextref\HUY
\nextref\MCK
\nextref\MCKA
\nextref\MI
\nextref\RU
\nextref\SAL
\nextref\SIK

\topmatter
\title 
An elliptic second main theorem
 \endtitle
\author  Julien Duval\footnote""{Laboratoire de math\'ematiques, Univ. Paris-Sud, Universit\'e Paris-Saclay, Orsay, France \newline julien.duval\@math.u-psud.fr\newline}
\footnote""{keywords : value distribution, J-holomorphic curve, elliptic structure \newline AMSC : 32H30, 32Q65.\newline}
\endauthor
  \endtopmatter 
\document
 \subhead 0. Introduction \endsubhead

\null\noindent  
Let $M$ be the projective plane $P^2(\C)$ endowed with an almost complex structure $J$ tamed by the Fubini-Study metric. Since Gromov [\GRO] we know that $M$ has plenty of $J$-lines (spheres of degree $1$ tangent to $J$). They build a projective plane $M'$, the dual of $M$, and satisfy the same incidence relations as the usual lines in $P^2(\C)$. 

\null\noindent
In particular we may consider a configuration of 4 $J$-lines in general position. According to B. Saleur [\SAL] (building on [\DUV]) such a configuration satisfies Borel's theorem. In other words, any $J$-entire curve avoiding the configuration has to be linear (contained in a $J$-line). Here a $J$-entire curve is a $J$-holomorphic map $f:\C\to M$. Our motivation is to quantify this fact. We want to estimate the intersections of a nonlinear $J$-entire curve with the configuration. 

\null\noindent
In the complex case this is nothing but Cartan's second main theorem (see [\RU]). Here we aim to use Ahlfors's more geometric approach based on duality and the use of singular metrics (see again [\RU]). However J.-C. Sikorav [\SIK] conjectured and B. McKay [\MCK] proved that $M'$ never carries a natural almost complex structure (unless $J$ is integrable) but only an elliptic one, a kind of nonlinear almost complex structure. Therefore we better try and work from the very beginning in this more general setting.

\null\noindent
So let now $M$ be an elliptic projective plane, i.e. the projective plane $P^2(\C)$ endowed with an elliptic structure $E$ tamed by a (normalized) symplectic form $\omega$. By [\SIK] $M$ still has plenty of $E$-lines. They build another elliptic projective plane $M'$, the dual of $M$, and satisfy the usual incidence relations. 
Alternatively elliptic projective planes can be seen purely in terms of incidence geometry as the smooth projective planes which are regular. See [\MCKA] for this line of thought.

\null\noindent
Fix a configuration of $q$ such $E$-lines $l_i$ in general position. Let $f:\C\to M$ be a nonlinear $E$-entire curve of infinite area. Denote by $T$ its characteristic function (measuring the growth) and by $N_i^{[2]}$ its counting function with respect to $l_i$ (measuring the intersections with $l_i$). We have $T(r)=\int_1^r \int_{D_t}f^*\omega \ {dt\over t} $ and $N_i^{[2]}(r)= \int_1^r n_i^{[2]}(t) {dt\over t}$. Here $D_t$ denotes the disc of radius $t$ centered at $0$ in $\C$ and $n_i^{[2]}(t)   $ the intersection between $f(D_t)$ and $l_i$ with multiplicities truncated at level 2, meaning that each intersection counts once if transversal and twice otherwise. As in the complex case our result reads as follows.

\th {Second main theorem} We have $(q-3)T(r)\leq \sum N_i^{[2]}(r)+o(T(r))\ \vert$.
\fth

\noindent
Here $\vert$ means that the estimate holds outside a set of finite Lebesgue measure. In particular if $q\geq 4$ we get plenty of intersections between $f(\C)$ and the configuration.

\null\noindent This paper is organized as follows. We begin with some background on elliptic projective planes, value distribution and Poincar\'e-Lelong formulas. We next enter elliptic value distribution theory, comparing the growths of an entire curve and its dual, and proving a first main theorem. We then turn to the proof of the second main theorem following the method of Ahlfors. Finally an appendix makes explicit some of the almost complex structures involved. 

\null\noindent
Thanks to the referee for his suggestions.

\noindent
All objects or maps are smooth except otherwise mentioned.
 
\null
 
\subhead 1. Elliptic projective planes \endsubhead

\null\noindent
We follow [\SIK] for this section. Again an elliptic projective plane $M$ is the projective plane $P^2(\C)$ equipped with an elliptic structure $E$ tamed by a (normalized) symplectic structure $\omega$.

\null\noindent
{\bf a) Elliptic structures.} They are nonlinear almost complex structures given by their complex tangent lines. To be precise, an elliptic structure tamed by $\omega$ is a subbundle $E$ (with $2$-dimensional compact fibers) of the Grassmanian $G(TM)$ of oriented $2$-planes in $TM$, such that $\omega$ is positive on its planes. We also require the following condition on its fiber $E_m$ over a point $m$ in $M$. Let $P$ be a plane of the Grassmanian $G(T_mM)$. Then $T_PG(T_mM)=\text{Hom}(P,T_mM/P)$. We want
$$T_P(E_m)\setminus \{0\}\subset \text{Iso}^+(P,T_mM/P). $$ 
This condition forces $E_m$ to be a sphere and its planes $P$ to foliate $T_mM\setminus \{0\}$. Moreover we get natural complex structures on $P$ and $T_mM/P$ such that the condition reads as
$$T_P(E_m) =\text{Hom}_\C(P,T_mM/P). $$ 
Hence $E_m$ itself inherits a natural structure of Riemann surface. The complex structures $j_P$ of the planes $P$ fit together in a selfmap $J_m$ of $ T_mM$ such that $J_m^2= -Id$. This "almost complex structure" is no longer linear in general. It is smooth outside the origin but only Lipschitz at the origin.  

\null\noindent
{\bf b) $E$-curves.} An embedded $E$-curve in $M$ is a surface tangent to $E$. It inherits a natural structure of Riemann surface. More generally a (parametrized) $E$-curve is a map $f:\Sigma \to M$ such that $J\circ df=df\circ i$. Here $(\Sigma,i)$ is a Riemann surface and $J$ is the (nonlinear) almost complex structure associated to $E$. When $\Sigma=\C$ we speak of $E$-entire curves.

\null\noindent
Geometric properties of complex curves extend to $E$-curves. For instance (nonconstant) $E$-curves are immersions except on a discrete set. At a singular point they still have a tangent. Distinct $E$-curves intersect positively.

\null\noindent 
{\bf c) Duality.} An $E$-line is an $E$-sphere of degree $1$. As said $M$ has plenty of $E$-lines. They build a projective plane $M'$, the dual of $M$, and satisfy the same incidence relations as in the complex case. In particular the fiber of the elliptic structure can be seen as $E_m =\{T_m\ l\ \vert \ l\ \text{$E$-line through}\  m\}$. So we may identify $E$ with the incidence variety $I=\{(m,l)\in M\times M'\ \vert \ m\in l\}$ together with its first projection $p$. Here we see indifferently $l$ either as a line in $M$ or a point in $M'$.

\null\noindent
Conversely the $E$-lines through $m$ form a sphere of degree $1$ in $M'$. We still denote it by $m$. Define $E'\subset G(TM')$ by its fiber $E'_l=\{T_l\ m\ \vert \ m\in l\}$. It turns out that $E'$ is an elliptic structure for $M'$. Again we may identify it with $I$ together with its second projection $p'$. Moreover this structure is tamed by the form $\omega'=\int_M [m] d\mu(m)$ given by Crofton's formula. Here $[m]$ denotes the current of integration over the line $m$ in $M'$ and $d\mu$ is a smooth measure of mass $1$ on $M$. Hence $M'$ also is an elliptic projective plane.

\null\noindent
Let $f: \C \to M$ be a (nonconstant) $E$-entire curve. We know that at each point this curve has a tangent direction so a well defined tangent $E$-line. We get in this way its dual curve $f': \C \to M'$. It is itself an $E'$-curve (see [\SIK]). 

\null\noindent
{\bf d) Pl\"ucker relations.} Denote by $V,V'\subset TI$ the vertical distributions, the respective kernels of $dp,dp'$. By a) their fibers have a natural complex structure. Note that $V'$ can also be seen as the tautological line bundle of $E$ (and similarly $V$ for $E'$). Indeed we have $d_{(m,l)}p (V')=T_m�\ l$.  We want now to compute the Chern classes of $V,V'$. As the whole situation deforms back to the complex case (see [\SIK]) the computation reduces to this case. So we may take $M=P^2(\C)$ with its usual complex structure and $E=P(TM)$. The relative Euler sequence (see [\HUY]) reads now as
$$ 0 \to V'\to p^*TM \to V\otimes V'\to 0 $$

\noindent
where $V'$ is seen as the tautological line bundle of $E$. We infer that $c_1(V)+2c_1(V')=p^*c_1(M)$. Similarly $c_1(V')+2c_1(V)=p'{^*}c_1(M')$. Recall that $c_1(M)= 3h$ where $h$ is the class of a line in $M$ (and similarly for $c_1(M')$). Solving we get the 

\th { Pl\"ucker relations} We have
$ c_1(V')= 2p^*h-p'{^*}h'$ and  $c_1(V)= 2p'{^*}h'-p^*h$. \fth

\null\noindent 
\subhead 2. Value distribution \endsubhead

\null\noindent
We collect some estimates from classical value distribution theory. We refer to [\RU] for background. We directly work on the complex line. 

\null\noindent
{\bf a) Characteristic functions.} Let $\alpha$ be a real 2-form on $\C$. Define its characteristic function by $T_\alpha(r)=\int_1^r\int_{D_t} \alpha \ {dt\over t}$. Here $D_t$ is the disc (centered at $0$) of radius $t$. If we write $\alpha=dd^cu$ then $T_\alpha(r)=2\pi \oint_ru + O(1)$ by integrating twice. Here $\oint_r$ denotes the mean on the circle of radius $r $ and $d^c=i(\overline \partial - \partial)$. 
In particular when $u$ is bounded from above we get $$T_{dd^cu}(r)\leq O(1). $$

\noindent
Note that $T_\alpha$ still makes sense if $\alpha$ is merely a measure. Hence we may allow mild singularities for $u$.
 
\null\noindent 
{\bf b) Ricci forms.} Suppose now $\alpha$ positive, i.e. $\alpha=\lambda dx \wedge dy$ with $\lambda>0$, and of infinite total mass. Hence $\log(r) =o(T_\alpha(r))$. Define the Ricci form of $\alpha$ by $Ric(\alpha)=dd^c\log \lambda$.
Then by Jensen's inequality and a calculus lemma (see [\RU]) we have $$T_{Ric(\alpha)}(r)\leq O(\log(rT_\alpha(r)))\ \vert \leq o(T_\alpha(r))\ \vert. $$
Here again $\vert$ means that the estimate holds outside a set of finite Lebesgue measure.

\null\noindent
{\bf c) Ahlfors lemma.} Our $T_\alpha(r)$ is the integrated area of concentric discs for the conformal metric given by $\alpha$. We may also consider the corresponding integrated length of concentric circles $L_\alpha(r)= 2\pi \int_1^r (\oint_t \sqrt \lambda) dt$.
Then the integrated Ahlfors lemma (see [\MI]) reads as $$L_\alpha(r)= o(T_\alpha(r))\ \vert.$$

\noindent
As before we may allow mild singularities (or zeros) for $\lambda$ in these estimates as long as the expressions involved make sense.

\null\noindent 
\subhead 3. Poincar\'e-Lelong formulas \endsubhead

\null\noindent We present a way to obtain Poincar\'e-Lelong formulas in an almost complex setting. We refer to [\BOT],[\HUY] for background on the geometry of vector bundles.

\null\noindent {\bf a) Angular forms.} Let $p:L\to X$ be a (smooth) hermitian line bundle over a compact manifold $X$. We identify $X$ with the $0$-section of $L$. Choose a hermitian connection $D$. Take a local unitary trivialization $L\vert_U\simeq U\times \C$, $v\mapsto (x,t)$ given by a unitary section $e$. Then $De= \omega e$ where $\omega$ is a purely imaginary 1-form on $U$, the connection form. Under a change of trivialization $e'=\lambda e$ we have $\omega'=\omega +{d\lambda\over \lambda}$ and $t'={t\over \lambda}$. Hence $d\omega$ is invariant and gives rise to a global 2-form  on $X$, the curvature of the connection. We write it as $-2\pi i c$ so that $c$ belongs to the first Chern class $c_1(L)$. Besides ${dt\over t}+\omega$ also is invariant. Its imaginary part divided by $2\pi$ gives rise to a global 1-form $\eta$ on $L\setminus X$, the angular form of the connection. We have by construction
$$ d\eta=[X]-p^*c   $$

\noindent where $[X]$ is the current of integration on $X$.

\null\noindent {\bf b) Poincar\'e-Lelong formula for $L$.}  
We assume now that the total space of $L$ carries an almost complex structure $J$. We suppose moreover that $J$ coincides with the given complex structure tangentially to the fibers and that the multiplications $\mu_r: L\to L, v\mapsto rv$ are $J$-holomorphic for $r$ in $\R$. We infer that
$$dt\circ J=idt+O(t) \tag 1$$

\noindent 
in a local unitary trivialization. So $\text{Re }({dt\over t}\circ J)+\text{Im }({dt\over t})$ remains bounded. Hence by our previous discussion we get the following Poincar\'e-Lelong formula
$$ {1\over 4\pi}d^c_J\log \Vert.\Vert^2= \eta + \beta. $$
 
\noindent 
Here $\Vert.\Vert$ denotes the norm of the hermitian metric, $d^c_Ju=-du\circ J$ for a function $u$, $\eta$ is an angular form on $L$ and $\beta$ a bounded 1-form on $L\setminus X$, vanishing tangentially to the fibers of $L$. We have $\beta=O( dp )$.

\null\noindent Note that by differentiating we get the usual form of the Poincar\'e-Lelong formula up to an extra term $d\beta$.

\null\noindent
{\bf c) Poincar\'e-Lelong formula for a hypersurface.} Let $(X,J)$ be a compact almost complex manifold and $H\subset X$ an almost complex hypersurface. We identify a neighborhood of $H$ in $X$ with a neighborhood $V$ of the $0$-section in its normal bundle $NH$. From the exact sequence
$$ 0\to TH\to TX\vert_H\to NH\to 0$$
we infer that $J$ induces a natural complex structure on the fibers of $NH$. We equip it with a hermitian metric and a hermitian connection. In a local unitary trivialization $(1)$ still holds near $H$. Take a cut-off function $\chi$ around $H$ with support in $V$. Define the function $u$ on $X$ by $u=\chi \Vert.\Vert^2+1-\chi$. 
Then again we get a Poincar\'e-Lelong formula which reads as
$$ {1\over 4\pi}d^c_J\log u= \alpha + \beta $$

\noindent 
where $\alpha =\chi\eta$ (with $\eta$ an angular form on $NH$) and $\beta$ is a bounded 1-form on $X\setminus H$. Moreover
$$ d\alpha=[H]-\tau $$

\noindent
where $\tau$ is a Thom form on $NH$ (see [\BOT]). In particular $\tau$ is Poincar\'e-dual to $H$.

\null
\subhead 4. Pl\"ucker estimates \endsubhead
 
\null\noindent
We go back to our context. From now on we fix a nonlinear $E$-entire curve (of infinite area) $f:\C\to M$ in an elliptic projective plane $M$. Its dual $f':\C\to M'$ also is a nonlinear $E'$-entire curve. Associated to them are their characteristic functions $T=T_{f^*\omega}$ and $T'=T_{f'{^*}\omega'}$. We want to compare them. This will be done by translating the Pl\"ucker relations at the level of characteristic functions via the Poincar\'e-Lelong formula. 

\null\noindent
For this denote by $N'$ the counting function of the zeros of $df$, $N'(r)=\int_1^rn'(t) {dt\over t}$ where $n'(t)$ counts the zeros (with multiplicities) of $df$ in $D_t$. Call $N''$ the corresponding function for $df'$. The estimates read as follows.

\th {Pl\"ucker estimates} We have $T'+N'\leq 2T+o(T+T') \ \vert$ and $T+N''\leq 2T'+o(T+T') \ \vert$.
\fth

\noindent In particular $T$ and $T'$ are comparable. We have $T=O(T')\ \vert$ and $T'=O(T)\ \vert$.

\null\noindent
We focus on the first estimate (the second follows by duality). To prove it we work at the level of the incidence variety $I$ and involve the line bundle $V'$ (see \S  1). We start by describing the natural (partial) almost complex structures on them and the corresponding Poincar\'e-Lelong formula.

\null\noindent
{\bf a) The structure $\tilde J$.}
Denote by $F$ the map $(f,f'):\C\to I$. It is tangent to the distribution $\Theta=V\oplus V'\subset TI$. Recall that $V=\ker dp$ is the vertical distribution for $p:I\to M$ (and similarly for $V'$) and that they both are complex line bundles. Hence we get a natural (partial) almost complex structure $\tilde J$ defined on $\Theta$. It is tamed by $\Omega=p^*\omega+p'{^*}\omega'$. By construction $F$ is a $\tilde J$-entire map, it satisfies $dF\circ i=\tilde J\circ dF$. Also note that $T_{F^*\Omega}=T+T'$. Actually $\tilde J$ extends as a genuine almost complex structure on $I$ (see [\MCK]) but we won't use this fact.

\null\noindent
{\bf b) The structure $J'$.}
We have more. Recall that $V'$ can also be seen as the tautological bundle of $E$. So we have a map $F'=(F, df({\partial \over \partial z})):\C\to V'$. It turns out that $V'$ carries a (partial nonlinear) almost complex structure $J'$ for which $F'$ is $J'$-holomorphic. This structure is constructed as follows. An (embedded) $E$-curve $C$ in $M$ lifts to a $\tilde J$-curve $\tilde C$ in $I$. Hence it inherits a natural structure of Riemann surface which makes $TC$ a complex manifold, as well as $V'\vert _{\tilde C}\simeq(p\vert_{\tilde C})^*(TC)$. Varying the curve these various complex structures fit together in $J'$. This structure is only defined on $\Pi=(dq')^{-1}(\Theta)$ (except on $(dq')^{-1}(V)$). Here $q':V'\to I$ is the projection. Moreover $J'$ is nonlinear, yet linear over each $\tilde J$ -line ($\neq V$) in $\Theta$. See the appendix for an explicit description.

\null\noindent
{\bf c) The Poincar\'e-Lelong formula.} The scheme described in \S 3 b) still works in this more general context, with $\beta$ reflecting the features of $J'$. We get the following formula
$$ {1\over 4\pi}d_{J'}^c\log \Vert .\Vert^2=\eta\vert_{\Pi}+\beta.  \tag 2$$
Here $\Vert \ \Vert$ is the hermitian metric on $V'$ induced by $\omega$, $\eta$ an angular form on $V'$, and $\beta$ a nonlinear form defined only on $\Pi$ (except over $V$), yet linear over each $\tilde J$-line $\neq V$, and bounded in the sense that $\beta=O(dq')$ (see appendix for the boundedness).  Moreover $$d\eta=[I]-q'{^*}c'  $$
where $[I]$ is the current of integration on $I$ (the $0$-section of $V'$) and $c'$ is in $c_1(V')$.

\noindent
Note that $c'$ is cohomologous to $2p^*\omega-p'{^*}\omega'$ by the Pl\"ucker relations. So after modifying $\eta$ (and $\beta$) by a suitable form we may actually suppose that 

$$d\eta=[I]-q'{^*}(2p^*\omega-p'{^*}\omega').$$

\noindent
{\bf d) End of the argument.}
Pulling back $(2)$ by $F'$ and differentiating we get
$${1\over 4\pi }d d^c\log \Vert F'\Vert^2=F'{^*}[I]-2f^*\omega+f'{^*}\omega'+d F'{^*}\beta. $$

\noindent
Also note that $\Vert F'\Vert^2dx\wedge dy=f^*\omega$ by definition of the metric on $V'$. So we have
$$ {1\over 4\pi}Ric(f^*\omega)=F'{^*}[I]-2f^*\omega+f'{^*}\omega'+d F'{^*}\beta. \tag 3$$

\noindent
Now take the characteristic functions in this equality. The Ricci term is estimated by \S 2 b). We have $T_{Ric(f^*\omega)}\leq o(T)\ \vert$. For the $\beta$ term we use instead \S 2 c) since $T_{dF'{^*}\beta}=O(L_{F^*\Omega})$ by Stokes theorem and the boundedness of $\beta$. We get $T_{dF'{^*}\beta}= o(T_{F^*\Omega})=o(T+T')\ \vert$. Hence we end up with our first Pl\"ucker estimate.

\null  
\subhead 5. First main theorem \endsubhead

\null\noindent
We keep the context and notations of the previous section. Let also $l$ be an $E$-line in $M$. We want to compare $T$ and the counting function (with multiplicities) $N=T_{f^*[l]}$ relative to $l$. As in the classical case this involves an extra term, the proximity function $m$ to $l$, $m(r)=\oint_r \log {1\over d(f,l)}$ where $d$ is a distance on $M$ (of diameter 1). Our first main theorem reads as follows.

\th {First main theorem} We have $N+m=T+o(T)\ \vert$. \fth

\noindent
In particular $m$ is controlled by $T$. 

\null\noindent
To prove this we again work at the level of the incidence variety $I$ using the Poincar\'e-Lelong formula for the hypersurface $H=p^{-1}(l)$. 

\null\noindent
{\bf a) The Poincar\'e-Lelong formula.} We want to apply \S 3 c). Note that $H$ is $\tilde J$-holomorphic in the sense that $TH\cap \Theta $ is always a $\tilde J$-line, except on the lift $\tilde l$ of $l$ where actually $\Theta=TH$. So (except on $\tilde l$) we have an exact sequence

$$ 0\to TH\cap \Theta \to \Theta\vert_H\to NH\to 0.$$

\noindent
It induces a natural complex structure $j$ on the fibers of $NH$ except along $\tilde l$. It turns out that $j$ extends Lipschitz along $\tilde l$ (see appendix). In this more general context the scheme described in \S 3 c) still works. There exists a function 
$u: I\to [0,1]$ vanishing at the second order on $H$ such that
$$ {1\over 4\pi}d^c_{\tilde J}\log u= \alpha\vert_ \Theta + \beta. \tag 4$$

\noindent 
Here we have identified a neighborhood of $H$ in $I$ with a neighborhood of the $0$-section in $NH$, $\alpha=\chi\eta$ where $\chi$ is a cut-off function around $H$ and $\eta$ an angular form on $NH$, and $\beta$ is a bounded 1-form on $I\setminus H$, defined only on the distribution $\Theta$. To be precise $u$ is only Lipschitz where the regularity of $j$ drops, namely on $NH\vert_{\tilde l}$. Also $\alpha$ and $\beta$ lose regularity but remain bounded there. Moreover
$$ d\alpha=[H]-\tau $$

\noindent
where $\tau$ is Poincar\'e-dual to $H$, hence cohomologous to $p^*\omega$. Therefore after modifying $\alpha$ (and $\beta$) we may actually suppose that

$$ d\alpha=p^*([l]-\omega). $$

\noindent
{\bf b) End of the argument.} Pulling back $(4)$ by $F$ and differentiating we get

$$ {1\over 4\pi}dd^c\log u\circ F= f^*[l]-f^*\omega +dF^* \beta.  $$

\noindent
Now take the characteristic functions in this equality. As before we have $T_{dF^*\beta}=o(T+T')=o(T) \ \vert$ (since $T$ and $T'$ are comparable). Also note that ${1\over 2}\oint \log {1\over u\circ F}=m+O(1)$. Hence we end up with our first main theorem.

\null  
\subhead 6. Second main theorem \endsubhead

\null\noindent
We explain the proof of our second main theorem following the method of Ahlfors. It relies on the use of singular forms with poles along the configuration and requires the control of the growths of $f$ and $f'$ with respect to these forms.

\null\noindent We keep the context and notations of the previous two sections. Let now $(l_i)$ be a configuration of $q$ $E$-lines in general position in $M$. Associated to each of them we have a counting function $N_i$ (with multiplicities) and a proximity function $m_i$. Recall that $N'$ is the counting function (with multiplicities) of the zeros of $df$, and similarly $N''$ for $df'$. We will actually prove the following stronger statement. 

\th {Second main theorem (strong form)} We have \ $2N'+N''+(q-3)T\leq \sum N_i+o(T)\  \vert $. \fth

\noindent
Its proof splits into the following estimates 
$$ 2N'+(q-4)T+2T'\leq \sum (N_i+ m_i') + o(T) \ \vert \tag 5$$
$$ N''+T +\sum m_i'\leq 2T'+o(T) \ \vert .\tag 6$$

\noindent
Here $m_i'$ is the proximity function of $l_i$ understood now as a point in $M'$. So $m'_i(r)=\oint_r \log {1\over d'(f',l_i)}$ where $d'$ is a distance on $M'$ (of diameter 1).

\null\noindent
These estimates are obtained in the same way as the Pl\"ucker estimates by substituting $\omega$ and $\omega'$ with singular forms. We describe them now.

\null\noindent
{\bf a) The singular forms.} Recall that associated to each line $l_i$ we have a Poincar\'e-Lelong formula involving a nonnegative function $u_i$ on $I$ vanishing at the second order on $p^{-1}(l_i)$ (see \S 5). We may suppose $u_i$ small by scaling. We also consider a nonnegative function $v_i$ on $I$ vanishing at the second order on $p'{^{-1}}(l_i)$ where $l_i$ is now seen as a point in $M'$.

\null\noindent
As in the classical case define singular weights by $w=(\Pi {v_i\over u_i(\log u_i)^2})^{1\over  2}$ and $w'=\Pi{1\over v_i(\log v_i)^2}$. The singular forms are  $\hat \omega =w\ p^*\omega$ and $\hat \omega'=w'\ p'{^*}\omega'$. Put $\hat T=T_{F^*\hat \omega}$ and $\hat T'=T_{F^*\hat \omega'}$ the associated characteristic functions. It turns out that both of them are controlled by $T$ (see below). We have $\hat T=O(T)\ \vert $ and $\hat T'=O(T) \ \vert$.

\null\noindent
{\bf b) Proof of $(5)$.} We compute the Ricci form of $F^*\hat \omega$. We have
$Ric(F^*\hat \omega)=Ric(f^*\omega)+dd^c\log w\circ F$. Using $(3)$ we get
$${1\over 2\pi}Ric(F^*\hat \omega)=2F'{^*}[I]-4f^*\omega+2f'{^*}\omega'+2dF'{^*} \beta +{1 \over 2\pi}dd^c\log w\circ F$$
where $\beta$ is a (partial nonlinear) bounded form.

\noindent
Take the characteristic functions in this equality. Estimating the Ricci and $\beta$ terms as in \S 4 d) and taking into account the control of $\hat T$ by $T$ we get

$$ 2N'-4T+2T'+ \oint \log w\circ F \leq o(T) \ \vert. \tag 7$$

\noindent
It remains to make explicit the weight term. Recall that ${1\over 2} \oint \log {1\over u_i\circ F}=m_i+ O(1)$. Similarly ${1\over 2} \oint \log {1\over v_i\circ F}=m'_i+ O(1)$. Moreover $m_i=T-N_i+o(T) \ \vert$ by the first main theorem. Finally by Jensen's inequality $\oint \log\log {1\over u_i\circ F} \leq \log(\oint \log {1\over u_i\circ F}) \leq o(T) \  \vert$. Summing up we get $$ \oint\log w\circ F\geq qT-\sum (N_i+m'_i)  -o(T)\  \vert.$$

\noindent
Together with $(7)$ this gives $(5)$.

\null\noindent
{\bf c) Proof of $(6)$.} We now compute $Ric(F^*\hat \omega')$. We get
$${1\over 4\pi}Ric(F^*\hat \omega')=F''{^*}[I]-2f'{^*}\omega'+f{^*}\omega+dF''{^* \beta'} +{1 \over 4\pi}dd^c\log w'\circ F$$

\noindent
where $\beta'$ is a (partial nonlinear) bounded form.

\noindent
Taking the characteristic functions in this equality, estimating the Ricci and $\beta'$ terms as above and using the control of $\hat T'$ by $T$ we get  
$$ N''-2T'+T+ {1\over 2}\oint \log w'\circ F \leq o(T) \ \vert. \tag 8$$

\noindent Making explicit the weight we have ${1\over 2}\oint \log w'\circ F  \geq \sum m'_i-o(T)\ \vert$. Here again the $\log \log $ terms are estimated via Jensen's inequality and absorbed in $o(T)$. Together with $(8)$ this gives $(6)$.

\null\noindent
This finishes the proof of the theorem modulo the controls of $\hat T$ and $\hat T'$.

\null\noindent
{\bf d) Control of $\hat T'$.} It suffices to estimate the characteristic function of $F^*{p'{^*}\omega' \over v_i(\log v_i)^2}$ for each $i$. For short we drop the index. Its control relies on the following estimate
$$-dd^c\log (-\log (v\circ F)) \geq  \epsilon F^* {p'{^*}\omega'\over v(\log v)^2}-C F^*\Omega + dF^*\beta.   \tag 9$$

\noindent
Here $\beta$ is a bounded 1-form defined only on $\Theta$ and $\epsilon,C>0$ are respectively small and large.

\null\noindent
Indeed taking the characteristic functions in $(9)$ and using the results of \S 2 as before we get $$T_{F^*{p'{^*}\omega' \over v(\log v)^2}}= O(T_\Omega)=O(T+T')=O(T)\ \vert.$$

\noindent
We turn to the proof of $(9)$. We first make precise the construction of $v$. Write $S$ for the fiber $p'{^{-1}}(l)$. Note that $\Theta\vert_S=TS\oplus V\vert_S$. So  $V\vert_S$ can be seen as a subbundle of $NS$. Choose a metric on $NS$ which is hermitian for $\tilde J$ on $V\vert_S$ and take $v=\Vert.\Vert^2$ in a neighborhood of $S$ in $I$ (identified with a neighborhood of the $0$-section of $NS$). 

\null\noindent
Locally in an unitary trivialization of $NS$ we have $v=\vert z\vert ^2+\vert w\vert ^2$ where $(z,w)$ are the coordinates in the fiber. Write $i$ for the standard complex structure in the fiber. By construction $\tilde J=i$ on $V$ for $(v=0)$. So $d_{\tilde J}^c v=d^c v\vert_\Theta +O(v)$. 

\null\noindent
We want to estimate $-dd^c\log (-\log (v\circ F))=dF^*{d_{\tilde J}^c v\over v (-\log v)}$. 

\noindent
We work locally. By the preceding remark we have $${d_{\tilde J}^c v\over v(-\log v)}={d^c v\vert_\Theta\over v (-\log v)}+O(1).\tag 10$$

\noindent
Differentiating the right-hand side we get $${vdd^cv-dv\wedge d^cv\over v^2(-\log v)}+{dv\wedge d^cv\over v^2(\log v)^2}.$$ We estimate each term. 

\noindent
A computation shows that $vdd^cv-dv\wedge d^cv$ is a nonnegative form for $i$. We infer that $vdd^cv-dv\wedge d^cv\geq_{\tilde J}-O(v^{3\over 2})$. Here 
$\geq_{\tilde J}0$ means being a nonnegative form for $\tilde J$ on $\Theta$. 

\noindent
Also note that ${dv\wedge d^cv\over v}$ is a positive form for $i$ on $V$ near $(v=0)$. We infer that  
$dv\wedge d^cv\geq_{\tilde J} 2\epsilon v\ p'{^*}\omega'-O(v^{3\over 2})$ for $\epsilon>0$ small.

\null\noindent
Summing up we get
$$d({d^c v\over v(-\log v)})\geq_{\tilde J} 2\epsilon {\ p'{^*}\omega'\over v(\log v)^2} -O({1\over {v^{{1\over 2}}(-\log v)}})\geq_{\tilde J} \epsilon {\ p'{^*}\omega'\over v(\log v)^2} -O(1).$$

\noindent
Together with $(10)$ this gives $(9)$ after pulling back by $F$.

\null\noindent
{\bf e) Control of $\hat T$.} By the sum to product trick (see [\RU]) it suffices to estimate the characteristic function of $F^*{v_ip^*\omega \over u_i(\log u_i)^2}$ for each $i$. For short we drop the index. Its control follows from the estimate (same notations as in $(9)$)
$$-dd^c\log (-\log (u\circ F)) \geq  \epsilon\ F^* {v\ p^*\omega\over u(\log u)^2}-C F^*\Omega + dF^*\beta. \tag 11  $$

\noindent
To prove this we proceed as above. We investigate ${d_{\tilde J}^c u\over u(-\log u)}$ locally.
Recall that $u$ is a (small) nonnegative function on $I$ vanishing at the second order on $H=p^{-1}(l)$. Locally in an unitary trivialization of $NH$ we have $u=\vert z\vert ^2$ where $z$ is the coordinate in the fiber. Moreover $d_{\tilde J}^c u=d^c u\vert_\Theta +O(u)$. So 
 $${d_{\tilde J}^c u\over u (-\log u)}={d^c u\vert_\Theta\over u (-\log u)}+O(1).\tag 12$$

\noindent
Differentiating the right-hand side we simply get ${du\wedge d^cu\over u^2(\log u)^2}.$ We estimate it.

\noindent
For this recall that $\Theta$ becomes tangent to $H$ precisely along the lift $\tilde l=H\cap S$ and that $v$ vanishes at the second order on $S$. We infer that for $\epsilon>0$ small $du\wedge d^cu\geq_{\tilde J} 2\epsilon uv\ p^*\omega - O(v u^{3\over 2}).$

\noindent
So we get
$$d({d^c u\over u (-\log u)})\geq_{\tilde J} 2\epsilon {v\ p^*\omega\over u(\log u)^2}-O({v\over {u^{{1\over 2}}(-\log u)}})\geq_{\tilde J} \epsilon {v\ p^*\omega\over u(\log u)^2} -O(1).$$
\noindent
Together with $(12)$ this gives $(11)$ after pulling back by $F$.

\null\noindent
{\bf f) Initial statement.} The second main theorem as stated in the introduction follows from this strong form since 
 
$$\sum N_i-(2N'+N'')\leq \sum N_i^{[2]} $$

\noindent
which in turn is obtained by integrating the punctual inequality 

$$\sum \nu_i-(2\mu'+\mu'')\leq \sum \nu_i^{[2]}. \tag 13$$

\noindent
Here $\nu_i$ is the intersection number between the entire curve and $l_i$, 
$\nu_i^{[2]}$ the same truncated at level $2$, $\mu'$ the multiplicity of the zero of $df$ and $\mu''$ the same for $df'$, all at a given point $z$. The truncation is understood geometrically. If $f(z)$ hits $l_i$, then $\nu_i^{[2]}=1$ if $f'(z)\neq l_i$ and  $\nu_i^{[2]}=2$ otherwise.

\null\noindent
The multiplicities $\mu',\mu''$ can also be seen geometrically. We refer to [\SIK] for background on singularities of $E$-curves. Define the multiplicity $k\ (\geq 1)$ of the curve at $z$ as its intersection number with any $E$-line through $f(z)$ different from its tangent line $f'(z)$. Denote also by $l\ (>k)$ the intersection number of the curve with its tangent line. As in the complex case we have $\mu'=k-1$. Moreover $l-k$ is the multiplicity of the dual curve at $z$ so $\mu''=l-k-1$. Hence  $2\mu'+\mu''=k+l-3$.

\null\noindent
The proof of $(13)$ goes by inspection, distinguishing at $f(z)$ the different possibilities for the position of the curve with respect to the configuration of lines.

\null\noindent
a) The curve hits a single line, transversally (the line differs from $f'(z)$). Then $(13)$ reads as $k-(k+l-3)\leq 1$.

\noindent
b) The curve hits a single line, tangentially (the line coincides with $f'(z))$. Then $(13)$ reads as $l-(k+l-3)\leq 2$.

\noindent
c) The curve hits a double point, transversally to the 2 lines. Then $(13)$ reads as $(k+k)-(k+l-3)\leq 1+1$.

\noindent
d) The curve hits a double point, tangentially to 1 of the 2 lines. Then $(13)$ reads as $(k+l)-(k+l-3)\leq 1+2$.

\null
\subhead 7. Appendix \endsubhead

\null\noindent
We make precise some features of the almost complex structures $J'$ and $j$. We work locally. We first describe the elliptic structure $E$ in coordinates.

\null\noindent
{\bf a) The structure $E$.}  We refer to [\SIK] for this paragraph. Let $(z,w)$ be coordinates on $M$ and $(\delta z, \delta w)$ coordinates in the fiber of $TM$. Take coordinates $(\lambda,\mu)$ in the fiber of the Grassmanian $G(TM)$ so that planes transversal to $(\delta z =0)$ can be written $(\delta w=\lambda \delta z +\mu \overline{\delta z})$. We normalize the elliptic structure as follows. We suppose that $(w=0)$ is an $E$-curve and that the induced complex structure on it coincides with the standard complex structure $i$. We also assume that the induced complex structure on the normal bundle of $(w=0)$ coincides with $i$ on the verticals $(z=cst)$. Then near the lift of $(w=0)$ the structure $E$ ($\subset G(TM)$) is given by the equation ($\mu=h(z,w,\lambda)$) for some function $h$. In particular $(z,w,\lambda)$ are coordinates for $E$. Moreover our normalizations translate to $(h=h_\lambda=h_{\overline \lambda}=0) $ along $(w=\lambda=0)$ (see [\SIK]). Here we use the notation $h_\lambda$  for the partial derivative of $h$ with respect to $\lambda$.

\null\noindent
We make explicit now the complex structures induced by $E$ on its planes $P$ ($\subset T_mM$) and on $T_mM/P$. Let $P$ be the plane of coordinates $(z,w,\lambda)$. For short, write $\alpha=h_\lambda(z,w,\lambda) $ and $\beta =h_{\overline \lambda}(z,w,\lambda)$. Identify $\text{Hom}(P,T_mM/P)$ with $\text{Hom}(\C,\C)$ (the $\C's$ being the $\delta z$-axis and the $\delta w$-axis). Then $T_PE\setminus \{0\}$ coincides with the set of automorphisms $L_s:\C\to\C$, $u\mapsto s u+( \alpha s +  \beta \overline s )  \overline u$ (for $s\in \C^*$). As said (see \S 1 a) there are unique (positively oriented) complex structures on the $\C's$ for which all the $L_s$ become complex linear. A way to see them is as follows.
\th{Lemma } There are unique (orientation preserving) automorphisms $A,B$ of $\C$ of the form $A(u)=u + a \overline u$, $B(u)=u + b\overline u$ such that all the $B\circ L_s\circ A$ are $\C$-linear. The $a,b$ are given by the equations
$( a+\alpha+ab\   \overline\beta =0)$ and
$(b+ \beta +ab\  \overline \alpha=0).$
\fth

\noindent The proof is left to the reader. We get a function $a$ related to the complex structure induced by $E$ on its planes $P$, and a function $b$ related to the complex structure on the quotients $TM/P$. Because of our normalizations $a$ and $b$ vanish on $(w=\lambda=0)$. By the lemma we have $a=-h_\lambda+O(\lambda^3)$ and $b=-h_{\overline \lambda}+O(\lambda^3)$ on $(w=0)$.

\null\noindent
{\bf b) The structure $J'$.} 
We make explicit the almost complex structure $J'$ on the tautological bundle $L$ of $E$ (see \S 4 b). We keep the notations and normalizations of  the previous paragraph. As already said $J'$ is obtained on $L\vert_{\tilde C}$ via the complex structure of $TC$ induced by $E$. Here $C$ is an $E$-curve and $\tilde C$ its lift in $E$. 
By definition the fiber of $L$ at the point $(z,w,\lambda)$ is the plane $(\delta w=\lambda\delta z+h(z,w,\lambda)\overline {\delta z})$. It is parametrized by $\delta z$, so $(z,w,\lambda, t=\delta z)$ are coordinates for $L$, and so are $(\delta z, \delta w,\delta \lambda, \delta t$) in the fiber of $TL$. 

\null\noindent
 Let $C$ be an $E$-curve tangent to $(w=0)$ at the origin. Write $C$ as $(w=f(z))$ with $f(0)=f_z(0)=f_{\overline z} (0)=0$. We have $(f_{\overline z}=h(z,f,f_z))$ as $C$ is an $E$-curve. Differentiating we get $ f_{z\overline z}(0)=  f_{\overline z^2}(0)=0$. On the other hand we may prescribe the remaining term $ r= f _{z^2}(0)$ of the $2$-jet as we wish (see [\SIK]). 

\null
\null\noindent  
Now we reparametrize $C$ conformally with respect to its complex structure induced by $E$. Via the initial parametrization the latter gives a complex structure $j_C$ on the $z$-axis. Consider a local diffeomorphism $\phi$ of $\C$, tangent to the identity at the origin, such that $$d\phi\circ i=j_C\circ d\phi. \tag 14$$ Hence $(\phi,f\circ\phi)$ is a conformal parametrization of $C$ and so is $\Phi=(\phi,f\circ\phi,f_z\circ\phi )$ for $\tilde C$. Differentiating we get a parametrization of $L\vert_{\tilde C}$ by $\Phi'(z,v)=(\Phi(z), d_z\phi(v))$. By definition of $J'$ we have $ d\Phi'\circ i=J'\circ d\Phi'$. Computing at the origin $O=(0,0,0)$ of $E$ we get
$$ J'_{(O,t)}(\delta z, 0, r \delta z,\delta t)=i(\delta z, 0, r \delta z,\delta t-2(t\phi_{z\overline z}(0) +\overline t  \phi _{ \overline z^2 } (0))\overline{\delta z}). $$ 
 We express now the second derivatives of $\phi$ in terms of $h$ and $r$. By the lemma above $(14)$ translates to $(\phi_{\overline z}=(a\circ \Phi)  \overline{ \phi_z}) $
where $(a+ h_\lambda +ab\  \overline {h_{ \overline\lambda}}=0)$. Differentiating and taking into account our normalizations we get  $\phi_{z\overline z}(0)=-r \ h_{\lambda^2}(O)$ and $ \phi _{ \overline z^2 } (0)=-\overline r\ h_{\lambda \overline \lambda}(O)$. So we end up with
$$ J'_{(O,t)}(\delta z, 0, r \delta z,\delta t)=i(\delta z, 0, r \delta z,\delta t+2(rt\ h_{\lambda^2}(O)+\overline{rt}\ h_{\lambda \overline \lambda}(O) ) \overline{\delta z}). \tag 15$$  

\noindent This is what we expect. Indeed at the origin we have $\Theta=(\delta w=0)$ and $V=(\delta z=\delta w=0)$. So $J'$ is linear over each line $(\delta \lambda =r\delta z)$ but nonlinear over $\Theta$ and actually not defined over $V$. 

\null\noindent
We also get the control of the form $\beta$. It follows by estimating ${dt\circ J'-idt \over t}$ (see \S 3 b). By $(15)$ we have $(dt\circ J'-idt )(\delta z,0,r\delta z,\delta t)=2i(t\ h_{\lambda^2}(O)\ r\overline{\delta z}+\overline t\ h_{\lambda \overline \lambda}(O) \ \overline{r\delta z})$
at $(O,t)$. Hence
 $${\vert(dt\circ J'-idt )(\delta z,0,\delta \lambda,\delta t)\vert \over \vert t\vert}\leq 2(\vert h_{\lambda^2}(O)\vert + \vert  h_{\lambda \overline \lambda}(O)\vert )\ \vert \delta \lambda\vert. $$

\noindent It is a bounded form.

\null\noindent
{\bf c) The structure $j$.} 
We briefly mention why the complex structure $j$ on the fibers of $NH$ extends Lipschitz on $\tilde l$ (see \S 5 a). Here $l$ is an $E$-line in $M$, $\tilde l$ its lift in $E$ and $H=p^{-1}(l)$. The structure $j$ is obtained from $\tilde J$ by projecting $\Theta$ onto the normal bundle of $H$, except on $\tilde l$ where $\Theta$ becomes tangent to $H$. We work locally, keeping the notations and normalizations of \S7 a). 

\null\noindent We identify $l$ with $(w=0)$, so $\tilde l=(w=\lambda=0)$. Let $P$ (in $H$) be the plane of coordinates $(z,0,\lambda)$. The structure $j$ on the fiber $NH_P$ can be seen by projecting horizontally the structure $j_P$ induced by $E$ on $P$ on the $\delta w$-axis. By the lemma above we already know the structure obtained by projecting vertically $j_P$ on the $\delta z$-axis. It is induced from $i$ by $A(u)=u+a\overline u$ where $a$ is as in the lemma. As $P=(\delta w=\lambda \delta z+h\overline {\delta z})$ we infer that $j$ is induced from $i$ by $C(u)=u+c\overline u$ where $c={h+\lambda a\over \overline \lambda+\overline h  a}$. Using the normalizations we get that $c=O(\lambda)$ so $j$ can be extended continuously through $\tilde l$ by $i$. Moreover it can be checked that the first derivatives of $c$ remain bounded, meaning that this extension is Lipschitz.

\Refs  
\widestnumber\no{99}
\refno=0
\bref \by R. Bott and L. Tu \book Differential forms in algebraic topology \bookinfo Graduate Texts in Mathematics\vol82  \publ 
 Springer \yr 1982 \publaddr Berlin 
 \endref

\bref \by J. Duval \paper Un th\'eor\`eme de Green presque complexe \jour Ann. Inst. Fourier \vol54\yr2004\pages2357--2367
 \endref

\bref \by M. Gromov \paper Pseudo holomorphic curves in symplectic 
 manifolds \jour Invent. Math. \vol82\yr1985\pages307--347
 \endref

\bref \by D. Huybrechts \book Complex geometry. An introduction. \bookinfo 
 Universitext \publ 
 Springer \yr 2005 \publaddr Berlin 
 \endref

\bref \by B. McKay \paper Dual curves and pseudoholomorphic curves \jour Selecta Math. \vol9\yr2003\pages251--311
\endref

\bref \by B. McKay \paper Smooth projective planes \jour Geom. Dedicata \vol116\yr2005 \pages157--202
\endref

\bref \by J. Miles \paper A note on Ahlfors' theory of covering surfaces \jour Proc. Amer. Math. Soc. \vol21\yr1969\pages30--32
\endref

\bref \by M. Ru \book Nevanlinna theory and its relation to Diophantine approximation  \publ World Scientific  \yr 2001 \publaddr River Edge
\endref 

\bref \by B. Saleur\paper Un th\'eor\`eme de Bloch presque complexe \jour Ann. Inst. Fourier\vol64\yr2014\pages401--428
\endref

\bref \by J.-C. Sikorav \book Dual elliptic planes, {\rm in} Actes des journ\'ees math\'ematiques \`a la m\'emoire de Jean Leray
\pages185--207 \bookinfo S\'emin. Congr. \vol9 \publ SMF \yr2004 \publaddr Paris
\endref

\endRefs

\enddocument